\documentclass[twoside,12pt]{article}
 \usepackage{bbm}
 \usepackage{mathrsfs}
  \usepackage{amsfonts}
\usepackage{amsmath}
\newcommand{\bref}[5]{\bibitem{#1} {#2} {\it #3} {\bf #4}#5.}

\begin{document}
  \title{A note on the third invariant factor of the Laplacian matrix of a graph
 \thanks{Supported by NSF of the People's Republic of China(Grant
 No. 10871189 ).}
  }
 \author{Jian Wang,  Yong-Liang Pan\thanks{Corresponding author. Email: ylpan@ustc.edu.cn}\,\,
 \\
  {\small Department of Mathematics, University of Science and Technology
         of China}\\
  {\small Hefei, Auhui 230026, The People's Republic of China}\\
}
\date{}
\maketitle {\centerline{\bf\sc Abstract}\vskip 8pt Let $G$ be a
simple connected graph with $n\geq 5$ vertices. In this note, we
will prove that $s_3(G)\leq n$, and characterize the graphs which
satisfy that $s_3(G)=n,\, n-1,\, n-2,\,$ or  $n-3$, where $s_3(G)$
is the third invariant factor of the Laplacian matrix of $G$.

\par \vskip 0.5pt {\bf Keywords}  Graph; Laplacian matrix;    Invariant factor; Smith
normal form.

{\bf 1991 AMS subject classification:}   15A18, 05C50 \\

\indent   Let $G=(V, E)$ be a simple connected graph with vertex set
$V=V(G)=\{v_1,\cdots,v_n\}$ and edge set $E=E(G)$. Denote the degree
of vertex $v_i$ by $d_i$ and let $D(G)=$diag$(d_1,\cdots, d_n)$. The
Laplacian
matrix  is  $L(G)=D(G)-A(G)$, where $A(G)$ is the $(0,1)-$adjacency matrix of $G$.\\
\indent Denote by $\Delta_i(G)$ the $i-$th determinantal minors of
$L(G)$, i.e., the greatest common divisor of all the $i-$by$-i$
determinantal minors of $L(G)$. Of course
$\Delta_i(G)\mid\Delta_{i+1}(G)$, $0<i<n$. The invariant factors of
$L(G)$ are defined by
$s_{i+1}(G)=\frac{\Delta_{i+1}(G)}{\Delta_{i}(G)}$, $0\leq i< n$,
where $\Delta_0(G)=1$. It is easy to see that  $s_i(G)|s_{i+1}(G),\
1\leq i\leq n-1$, and $s_n(G)=0$ since $L(G)$ is singular. The Smith
normal form of $L(G)$ is  the $n-$ square diagonal matrix $F(G)$
whose $(i,i)$ entry is $s_i(G)$. It follows from the well known
matrix-tree theorem  that $\Delta_{n-1}(G)=s_1(G)s_2(G)\cdots
s_{n-1}(G)$ is equal to the spanning tree number of $G$. So  the
invariant factors of $G$ can be used to distinguish pairs of
non-isomorphic graphs which have the same  spanning tree number,
and so there is considerable interest in their properties.\\
\indent Since $G$ is a simple  graph, its invariant factor $s_1(G)$
must be equal to 1, however most of the others are not easy to be
determined. From the following lemma, we know that $s_2(G)=1$ if $G$
is not the complete graph $K_n$, while
$F(K_n)=$diag$(1,n,\cdots,n,0)$. In this note, we will show that
$s_3(G)\leq n$, and
characterize the graphs which satisfy that $s_3(G)=n,\, n-1,\, n-2,\,$ or  $n-3$.\\

\noindent{\bf Lemma }\, ([1]) For a simple connected graph $G$ with
order $n\geq 3$, $s_2(G)\not=1$ if and
only if $G$ is the complete graph $K_n$, which has $s_i(G)=n,\, 2\leq i\leq n-1$.\\

\indent  In the following theorem, $v\cdot G$ denotes the graph
obtained by adding an edge joining some  vertex of $G$ to a further
vertex $v$; $G-2e$ denotes the graph obtained from $G$ by deleting
two edges which have no common vertex; $G-C_4$ denotes the graph
obtained from $G$ by deleting a circle of length 4; $G-2C_3$ denotes
the graph obtained from $G$ by deleting 6 edges
in two cycles of length 3 which have no common vertices (See Fig. 1). In the proof of the following theorem,
$x \sim y $ means that the vertices $x$ and $y$ are
adjacent and $x\not\sim y$ means that they are
not adjacent.\\
\\
\setlength{\unitlength}{.75cm}
\begin{picture}(6,6)(0,0)
\put(1,0){$\circ$}\put(5,0){$\circ$}
\put(0,2){$\circ$}\put(6,2){$\circ$}\put(3,5){$\circ$}
\put(1,4){$\circ$}\put(5,4){$\circ$}
\put(1.15,0.25){\line(2,5){1.93}}\put(1.19,0.22){\line(1,1){3.85}}
\put(1.22,0.18){\line(5,2){4.85}}\put(1.25,0.1){\line(1,0){3.8}}
\put(.18,2.25){\line(1,1){2.88}}\put(.2,2.2){\line(5,2){4.85}}
\put(.25,2.12){\line(1,0){5.8}}\put(.2,2.05){\line(5,-2){4.85}}
\put(1.2,4.22){\line(2,1){1.85}}\put(1.25,4.15){\line(1,0){3.8}}
\put(1.25,4.05){\line(5,-2){4.8}}\put(1.15,4.02){\line(1,-1){3.9}}
\put(3.23,5.1){\line(2,-1){1.8}}\put(3.18,5.05){\line(1,-1){2.9}}
\put(3.15,5.03){\line(2,-5){1.92}} \put(2,-2 ){ $K_7-2C_3$}
\end{picture}
\hfill \setlength{\unitlength}{.75cm}
\begin{picture}(6,6)(3,0)
\put(2,0){$\circ$}\put(2,4){$\circ$}\put(6,0){$\circ$}\put(6,4){$\circ$}
\put(0,2){$\circ$}\put(4,2){$\circ$}\put(8,2){$\circ$}
\put(1.9,4.4){$v_1$}\put(5.9,4.4){$v_2$}\put(3.9,2.4){$v_3$}
\put(1.9,-.3){$v_4$}\put(5.9,-.3){$v_5$}
\put(-.5,2){$v_6$}\put(8.3,2){$v_7$}
\put(2.25,0.2){\line(1,1){1.85}}\put(2.12,0.25){\line(0,1){3.8}}
\put(6.25,0.2){\line(1,1){1.85}}\put(6.12,0.25){\line(0,1){3.8}}
\put(0.18,2.07){\line(1,-1){1.87}}\put(.2,2.25){\line(1,1){1.85}}\put(.25,2.12){\line(1,0){3.8}}
\put(4.18,2.07){\line(1,-1){1.87}}\put(4.2,2.25){\line(1,1){1.85}}\put(4.25,2.12){\line(1,0){3.8}}
\put(2.2,4.07){\line(1,-1){1.87}}\put(6.2,4.07){\line(1,-1){1.86}}
\put(2,-2){ $K_7-K_{3,3}$}
\end{picture}\\
\setlength{\unitlength}{1cm}
\begin{picture}(6,6)(0,0)
\put(1,1){$\circ$}\put(1,3){$\circ$}
\put(6,0){$\circ$}\put(6,2){$\circ$}\put(6,4){$\circ$}
\put(1.18,1.12){\line(5,1){4.85}}\put(1.18,1.08){\line(5,-1){4.85}}
\put(1.18,3.12){\line(5,1){4.85}}\put(1.18,3.08){\line(5,-1){4.85}}
\put(1.16,1.15){\line(5,3){4.85}}\put(1.15,3.05){\line(5,-3){4.89}}
\put(1.1,1.18){\line(0,1){1.85}} \put(2.5,-1){ $K_5-C_3$}
\end{picture}\hfill
\setlength{\unitlength}{1cm}
\begin{picture}(6,6)(1,0)
\put(1,0){$\circ$}\put(6,0){$\circ$}
\put(1,2){$\circ$}\put(6,2){$\circ$}
\put(1,4){$\circ$}\put(6,4){$\circ$}
\put(1.17,0.12){\line(5,2){4.85}}\put(1.15,0.15){\line(5,4){4.87}}\put(1.19,0.1){\line(1,0){4.83}}
\put(1.16,2.13){\line(5,2){4.85}}\put(1.16,2.07){\line(5,-2){4.87}}\put(1.19,2.09){\line(1,0){4.83}}
\put(1.18,4.08){\line(5,-2){4.85}}\put(1.17,4.03){\line(5,-4){4.87}}\put(1.19,4.1){\line(1,0){4.83}}
\put(1.1,0.18){\line(0,1){1.85}}\put(1.1,2.18){\line(0,1){1.85}}
\qbezier(1.05,.15)(0,2.15)(1.05,4.08) \put(2,-1){ $K_6-C_3$}
\end{picture}\\
\setlength{\unitlength}{1cm}
\begin{picture}(6,8)(-0.5,0)
\put(1,1){$\circ$}\put(1,5){$\circ$}
\put(5,1){$\circ$}\put(5,5){$\circ$} \put(3,3){$\circ$}
\put(1.16,1.16){\line(1,1){1.89}}\put(1.18,1.11){\line(1,0){3.85}}
\put(1.16,5.04){\line(1,-1){1.87}}\put(1.18,5.11){\line(1,0){3.85}}
\put(3.16,3.16){\line(1,1){1.89}}\put(3.16,3.04){\line(1,-1){1.87}}
\put(1.1,1.18){\line(0,1){3.85}}\put(5.1,1.18){\line(0,1){3.85}}
\put(2.5,0){ $K_5-2e$}
\end{picture}\hfill
\setlength{\unitlength}{1cm}
\begin{picture}(6,8)(0.5,0)
\put(1,1){$\circ$}\put(1,5){$\circ$}
\put(5,1){$\circ$}\put(5,5){$\circ$} \put(3,3){$\circ$}
\put(1.16,1.16){\line(1,1){1.89}}\put(1.18,1.11){\line(1,0){3.85}}
\put(1.16,5.04){\line(1,-1){1.87}}\put(1.18,5.11){\line(1,0){3.85}}
\put(3.17,3.17){\line(1,1){1.89}}\put(3.16,3.04){\line(1,-1){1.88}}
\put(2.5,0){ $K_5-C_4$} \put(-1.5,-1){Fig. 1}
\end{picture}
\vspace{1.5cm}\\
\noindent{\bf Theorem } Let $G\not=K_n$ be a simple connected graph
with order $n\geq 5$. Then $s_3(G)\leq n$. Moreover, $s_3(G)=n$ if
and only if $G=K_n-e$, where $e$ is an edge of $K_n$; $s_3(G)=n-1$
if and only if $G=v\cdot K_{n-1}$; $s_3(G)=n-2$ if and only if $n=5$
and $G=K_5-2e$ or $G=K_5-C_4$; $s_3(G)=n-3$ if and only if $G$ is
one of the following 6 graphs: $K_{2,3}$,
$K_5-C_3$, $K_6-C_3$,  $K_7-2C_3$, $K_{3,3}$ and $K_7-K_{3,3}$.\\

\noindent{\bf Proof} Since $G\not=K_n$, then its diameter is at
least 2. In fact, we only need to consider the graphs with diameter
2, since if the diameter of $G$ is more than 2 then by theorem 4.5
in [1] we have that $s_3(G)=1$. Let $v_1$ and $v_2$ be two
nonadjacent vertices in $G$. There is a further
vertex $v_3$ which is adjacent to both $v_1$ and $v_2$. Now we need to distinguish 8 cases to go on the argument.\\
\indent $Case$ 1.\;  Some vertex  $v_4$ in $V(G)/\{v_1, v_2, v_3\}$
satisfies that $v_4\not\sim v_1$, $v_4\not\sim v_2$, $v_4\not\sim
v_3$. Since $G$ is connected,  there is some  vertex $v_5$ in
$V(G)/\{v_1, v_2, v_3, v_4\}$  adjacent to both $v_1$ and $v_4$.
Clearly,  $\det (L[1,3,4|2,3,5])=-1$, where $L[1,3,4|2,3,5]$ is the
submatrix of $L(G)$ that lies in the rows corresponding to vertices
$v_1, v_3, v_4$
and columns corresponding to vertices $v_2, v_3, v_5$. Therefore $s_3=1$.\\
\indent $Case$ 2.\; Some vertex  $v_4$ in $V(G)/\{v_1, v_2, v_3\}$
satisfies that $v_4\not\sim v_1$, $v_4\not\sim v_2$, $v_4\sim v_3$.
In this case, $\det(L[1, 2, 3|2, 3, 4])=d_2\leq n-2$. So  $s_3\leq n-2$.\\
\indent $Case$ 3.\;  Some vertex  $v_4$ in $V(G)/\{v_1, v_2, v_3\}$
satisfies that  $v_4\not\sim v_1$, $v_4\sim v_2$, $v_4\not\sim v_3$.
In this case,  $\det(L[1, 2, 3|2, 3, 4])=1$, and hence $s_3=1$.\\
\indent $Case$ 4.\; Some vertex  $v_4$ in $V(G)/\{v_1, v_2, v_3\}$
Satisfies that $v_4\not\sim v_1$, $v_4\sim v_2$, $v_4\sim v_3$.
In this case,  $|\det(L[1, 2, 3|2, 3, 4])|=d_2+1\leq n-1$. Hence $s_3\leq n-1$.\\
\indent $Case$ 5.\; Some vertex  $v_4$ in $V(G)/\{v_1, v_2, v_3\}$
satisfies that $v_4\sim v_1$, $v_4\not\sim v_2$, $v_4\not\sim v_3$.
In this case, very similar to $case$ 3, we have $s_3=1$.\\
\indent $Case$ 6.\; Some vertex  $v_4$ in $V(G)/\{v_1, v_2, v_3\}$
satisfies that $v_4\sim v_1$, $v_4\not\sim v_2$, $v_4\sim v_3$.
In this case,  very similar to $case$ 4, we have $s_3\leq n-1$.\\
\indent $Case\, 1 -Case\, 6$ show that if some vertex $v_4$ in
$V(G)/\{v_1, v_2, v_3\}$ is not adjacent to both $v_1$ and $v_2$,
then $s_3<n$. So, we will only need to deal with  the cases in which
every further vertex in
$V(G)/\{v_1, v_2, v_3\}$ is adjacent to both  $v_1$ and $v_2$.\\
\indent $Case$ 7.\; Every vertex  in $V(G)/\{v_1, v_2, v_3\}$ is
adjacent to both $v_1$ and $v_2$, and  at least one vertex $v_4$ is
not adjacent to
$v_3$. In this case, we distinguish  3 subcases.\\
\indent $Subcase$ 1.\, There is some vertex $v_5$ in
$V(G)/\{v_1,v_2,v_3,v_4\}$  adjacent to all of the vertices $v_1,
v_2$ and $v_3$. Then
$\det(L[1, 2, 3|1, 4, 5])=d_1=n-2$. Hence $s_3\leq n-2$.\\
\indent $Subcase$ 2.\, Every vertex in $V(G)/\{v_1,v_2,v_3\}$ is
not adjacent to $v_3$, and the induced subgraph $G[v_4, \cdots,
v_n]\neq K_{n-3}$. If we choose any two nonadjacent vertices in
$\{v_4, \cdots, v_n\}$  as $v_4$ and $v_5$, then we have
$-\det(L[1,4,5|1,3,5])=-\begin{vmatrix}
d_1&-1& -1\\
-1&0&0\\
-1&0&d_5\end{vmatrix}=d_5$.
Hence we have that $s_3\leq d_5\leq n-3$.\\
\indent $Subcase$ 3.\,  Every vertex in $V(G)/\{v_1,v_2,v_3\}$ is
not adjacent to $v_3$, but $G[v_4, \cdots, v_n]\\=K_{n-3}$.
It is not difficult to obtain that  $F(G)=$diag$(1, 1, 1, n-1, \cdots, n-1,\, 2(n-1)(n-2), 0)$.\\
\indent $Case$ 8.\; Every vertex   $G-\{v_1, v_2, v_3\}$ is adjacent
to all of the vertices $v_1$,
$v_2$ and $v_3$. In this case, we distinguish  two subcases.\\
\indent $Subcase$ 1.\, $G-\{v_1,v_2,v_3\}\neq K_{n-3}$, then there
are two nonadjacent vertices $v_4$ and $v_5$ in $V(G)/\{v_1, v_2,
v_3\}$. It follows that $\det(L[2,3,4|1,4,5])=d_4\leq n-2$. Hence
$s_3\leq n-2$. (In fact, if we regard the vertices $v_4$ as $v_1$,
$v_5$ as $v_2$, $v_1$ as $v_3$, $v_2$ as $v_4$, and $v_3$ as $v_5$,
then we are back in the subcase 1 of case 7.)\\
\indent $Subcase$ 2.\, $G[v_4, \cdots, v_n]=K_{n-3}$. Then
$G=K_n-e$. It is not difficult to obtain
$F(K_n-e)$=diag$(1, 1,  n, \cdots, n,\, n(n-2), 0)$.\\

\indent From above argument we have  that $s_3\leq n$ and $s_3=n$ if and only if $G$ is $K_n-e$.\\

\indent Clearly,  case 6 is symmetric to case 4,
the required graphs in case 4 are the isomorphic to the required graphs in case 4.\\

\indent From proposition 1 in [2], we know that  $F(v\cdot
K_{n-1})=$diag$(1, 1, n-1, \cdots, n-1, 0)$.
 Now  we prove the converse: if $s_3(G)=n-1$ then $G=v\cdot K_{n-1}$.\\
\indent  From the  argument of the above 8 cases, it follows  that
if $s_3(G)=n-1$ then every  vertex in $V(G)/\{v_1, v_2, v_3\}$ is
adjacent to $v_3$ and only $case$ 4 or $case$ 6 may occur. If $case$
4 occurs, then $\left|\det L [1,2,3|2,3,4]\right| =d_2+1\leq n-1$.
It follows from $s_3(G)=n-1$  that $d_2=n-2$ and then $case$ 6 will
never occur. Similarly,  if $case$ 6 occurs then  $case$ 4 will
never occur. Without loss of generality,
we assume that only  $case$ 4 occurs. We need to deal with two subcases here.\\
\indent $Subcase$ 1.\, There are two vertices in $\{v_4,\cdots,
v_n\}$ which  are not adjacent. We regard the two nonadjacent
vertices as $v_1$, $v_2$, and regard $v_1$ as $v_4$,
we are then back in $case$ 2, so we have that $s_3\leq n-2$.\\
\indent $Subcase$ 2.\,  $G[v_4, \cdots, v_n]=K_{n-3}$. Note that
$v_2$ and $v_3$ are adjacent to every vertex in
$V(G)/\{v_1, v_2, v_3\}$, thus $G=v\cdot K_{n-1}$.\\

\indent A direct calculation can show that $F(K_5-2e)=$diag $(1, 1,
3, 15, 0)$ and $F(K_5-C_4)=$diag$(1, 1, 3, 3, 0)$.
Now we prove that if $s_3=n-2$, then $n=5$ and $G=K_5-2e$ or $K_5-C_4$.\\
\indent By the above argument,
we know that if $s_3=n-2$ then $case$s 1, 3 and 5 may not occur. \\
\indent If $case$ 2 occurs, then
$\det(L[1, 2, 3|2, 3, 4])=d_2\leq n-2$. So $d_2$ must be $n-2$. It is a contradiction to  $case$ 2.\\
\indent If $case$ 4 occurs, then $\det(L[1, 2, 3|2, 3, 4])=d_2+1\leq
n-1$, and hence $d_2=n-3$. Then we must have exact one vertex $v_5$
in $case$ 6. If $n>5$, there is another vertex $v_6$ in $case$ 7 or
$case$ 8. We have, if $v_6\not \sim v_3$, then $\det(L[1, 2, 3|4, 5,
6])=2<n-2$; and if $v_6\sim v_3$, then $\det(L[1, 2, 3|4, 5,
6])=1<n-2$. So $n=5$.  Now,  $v_5\sim v_1$, $v_5\not\sim v_2$,
$v_5\sim v_3$. If $v_4\not\sim v_5$, then $G=K_5-C_4$, whose Smith
normal form is diag$\{1,1,3,3,0\}$; if $v_4\sim v_5$,
then $G=K_5-P_4$, whose Smith normal form is diag$\{1,1,1,21,0\}$. Impossible.\\
\indent If $case$ 7 occurs, then from the above argument, we know
that only its $subcase$ 1 occurs. Note that $\det(L[1, 2, 3|1, 4,
5])=\left|
  \begin{array}{ccc}
    d_1 & -1 & -1 \\
    0 & -1 & -1 \\
    -1 & 0 & -1 \\
  \end{array}
\right|=d_1\leq n-2$ and $\det(L[1, 2, 3|2, 4, 5])=\left|
  \begin{array}{ccc}
    0 & -1 & -1 \\
    d_2 & -1 & -1 \\
    -1 & 0 & -1 \\
  \end{array}
\right|=-d_2$, so $d_1=d_2=n-2$. Moreover,
 $\det(L[2, 3, 4|1, 3, 5])=\left|
  \begin{array}{ccc}
    0 & -1 & -1 \\
    -1 & d_3 & -1 \\
    -1 & 0 & x \\
  \end{array}
\right|=-(d_3+1+x)$ and $\det(L[2, 3, 4|1, 4, 5])=\left|
  \begin{array}{ccc}
    0 & -1 & -1 \\
    -1 & 0 & -1 \\
    -1 & d_4 & x \\
  \end{array}
\right|=d_4-1-x\leq n-2$, where $x$ is 0 if $v_4\not\sim v_5$, or
$-1$ if $v_4\sim v_5$. If $x=0$, then $d_4-1=n-2$ and it follows
that $d_4=n-1$, impossible. So $x=-1$, and then $d_3=d_4=n-2$. For
$i\geq 5$, $\det(L[1,4,i|2,3,i])=\left|
  \begin{array}{ccc}
    0 & -1 & -1 \\
    -1 & 0 & -1 \\
    -1 & -1 & d_i \\
  \end{array}
\right|=-(d_i+2)$. So $n-2\leq d_5+2\leq n+1$.\\
\indent $\bullet$\; If $d_5+2=n+1$, then $n-2$ divides $n+1$, thus $n=5$ and hence $G=K_5-2e$.\\
\indent $\bullet$\; If $d_5+2=n$, then $n-2$ divides $n$. Thus $n=4$, a contradiction.\\
\indent $\bullet$\; If $d_5+2=n-2$, then there are further 3 vertices $v_6$, $v_7$ and $v_8$  not
adjacent to $v_5$. Now  $\det(L[3, 5, 6|1,4,5])=\left|
  \begin{array}{ccc}
    -1 & 0 & -1 \\
    -1 & -1 & d_5 \\
    -1 & -1 & 0 \\
  \end{array}
\right|=-d_5=n-4$. Now $n-4$ divides $n-2$, it follows that $n=5$, or 6. Impossible.\\
\indent If $case$ 8 occurs,  then its $subcase$ 1 occurs and
$subcase$ 2 does not.
Then we only need to deal with $subcase$ 1 of $case$ 7, it has been done.\\

\indent With the aid of  Maple, we obtain the Smith normal forms of the graphs
$K_{2,3}$, $K_5-C_3$, $K_6-C_3$,  $K_7-2C_3$, $K_{3,3}$ and
$K_7-K_{3,3}$ as follows: $F(K_{2,3})=$diag$(1,1,2,6,0)$,
$F(K_5-C_3)=$diag$(1,1,2,10,0)$, $F(K_6-C_3)=$diag$(1,1,3,6,18,0)$,
$F(K_7-2C_3)=$diag\\$(1,1,4,4,4,28,0)$,
$F(K_{3,3})=$diag$(1,1,3,3,9,0)$, $F(K_7-K_{3,3})=$diag$(1,\ 1,\ 4,\
4,\ 4,\ 4,\ 0)$.
In the following, we will prove that if $s_3=n-3$  then $G$ must be one of these 6 graphs.\\
\indent By the above argument, we know that $cases$ 1, 3, 5 can not occur.\\
\indent If $case$ 2 occurs, then $\det(L[1, 2, 3|2, 3, 4])=d_2$,
$\det(L[1, 2, 3|1, 3, 4])=d_1$ and  $\det(L[2, 3, 4|1, 3, 4])=-d_4$.
Hence $d_1=d_2=d_4=n-3$.
Consider the number of vertices with degree $n-1$, we distinguish 3 subcases.\\
\indent  $Subcase$ 1.  $G$ has at least 3 vertices with degree
$n-1$, then $L(G)$ has a submatrix $L_1=\left(\begin{array}{cc}
(n-3)I_3 & -J_3\\
-J_3& nI_3-J_3
\end{array}\right)$, where $I_3$ is the $3\times 3$ identity matrix, $J_3$ is the $3\times 3$ all 1's matrix.
Note that $\det(L_1[1,4,6|2,4,5])= \left|
  \begin{array}{ccc}
    0 & -1 & -1 \\
    -1 & n-1 & -1 \\
    -1 & -1 & -1 \\
  \end{array}
\right|=-n.$ So $n-3$ divides $n$, it follows that $n=6$ and hence $G=K_6-C_3$.\\
\indent $Subcase$ 2.  $G$ has 1, or 2 vertices with degree $n-1$.\\
\indent $\bullet$ If $n=5$, then clearly, $G=K_5-C_3$.\\
\indent $\bullet$ If $n\geq 6$, then suppose $v_i\not\sim v_j$,
where $v_i$, $v_j\in V(G)/\{v_1,v_2,v_4\}$. $L(G)$ has a submatrix
$L_2=\left(\begin{array}{cc}
(n-3)I_3 & -J_3\\
-J_3& B
\end{array}\right)$, where $B=\left(
  \begin{array}{ccc}
     d_i & 0 & -1 \\
     0 & d_j & -1 \\
    -1 & -1 & n-1
  \end{array}\right).$
Then $|\det(L_2[1,4,6|2,4,5])|= |\left|
  \begin{array}{ccc}
    0 & -1 & -1 \\
    -1 & d_i & 0 \\
    -1 & -1 & -1 \\
  \end{array}
\right||=d_i\leq n-2$. So $d_i=n-3$. In the same way, we can get
$d_j=n-3$. Thus the  vertices of $G$ share two degrees: $n-1$ or
$n-3$. $\det(L_2[2,4,6|3,5,6])=\left|
  \begin{array}{ccc}
    0 & -1 & -1 \\
    -1 & 0 & -1 \\
    -1 & -1 & n-1 \\
  \end{array}
\right|
=-(n+1)$. Hence
$n-3$ divides $n+1$, then $n=7$ and it follows that $G=K_7-2C_3$.\\
\indent  $Subcase$ 3. $G$ has no vertex  with degree $n-1$.\\
\indent $\bullet$ If $n=5$, clearly $G=K_5-C_3-e=K_{2,3}$.\\
\indent $\bullet$ If $n=6$, clearly $G=K_6-2C_3=K_{3,3}$.\\
\indent $\bullet$ If $n\geq7$, then $L(G)$ has a principal submatrix
$L_3=\begin{pmatrix}(n-3)I_3 & -J_{3\times 4}\\
-J_{4\times 3}& C\end{pmatrix}$, where $C=\left(
  \begin{array}{cccc}
     d_i & 0 & y_1 & y_3 \\
     0 & d_j & y_2 & y_4 \\
     y_1 & y_2 & d_u & y_5 \\
     y_3 & y_4 & y_5 & d_v
  \end{array}\right)$
with $y_i=0$ or $-1$. Note that $\det(L_3[1, 4, 7|3, 5,
6])=y_1+y_4-y_5$, $\det(L_3[2, 4, 6|3, 5, 7])=y_2+y_3-y_5$. Now
$(n-3)\mid(y_1+y_4-y_5)$ and $(n-3)\mid(y_2+y_3-y_5)$, it follows
that  $y_1=y_2=y_3=y_4=y_5=0$ and hence $d_i\leq n-4$. Now
$\det(L_4[1, 4, 7|3, 4, 6])= \left|
  \begin{array}{ccc}
    0 & -1 & -1 \\
    -1 & d_i & 0 \\
    -1 & 0 & 0 \\
  \end{array}
\right|=-d_i$.  Therefore $n-3$ divides $d_i$. But  $d_i\leq n-4$, so it is impossible.\\
\indent If $case$ 4 occurs, then $|\det(L[1, 2, 3|2, 3,
4])|=d_2+1$. Therefore we have $n-3\leq d_2+1\leq n-1$. If $d_2+1=n-1$, then $n-3$ divides $n-1$, thus $n=5$
and $d_2=n-2=3$. So $v_5$ must be in case 4. If $v_4\sim v_5$, then
$G=v\cdot K_4$, whose $F(G)$=diag$(1,1,4,4,0)$. It is impossible.
If $v_4\not\sim v_5$, then a direct calculation can show $F(G)$=diag$(1,1,1,8,0)$,
it is a contradiction.
So $d_2+1=n-3$. There must be some vertex $v_5$ in  $case$ 6 and hence we have
$-\det(L[1, 2, 3|1, 3, 5])=-\left|
  \begin{array}{ccc}
    d_1 & -1 & -1 \\
    0 & -1 & 0 \\
    -1 & d_3 & -1 \\
  \end{array}
\right|=d_1+1\leq n-1$, so $d_1=n-4$. Then there are two vertices
$v_6$ and $v_7$  such that $v_7$ together with $v_5$ are in  $case$
6, and $v_6$ together with $v_4$ are in $case$ 4. Then
$-\det(L[1,2,3|3,6,7])=- \left|
  \begin{array}{ccc}
    -1  &  0 & -1 \\
    -1  & -1 & 0 \\
    d_3 & -1 & -1 \\
  \end{array}
\right|=d_3+2.$ Thus we have that $n-3$ divides $d_3+2$ and $n-3\leq d_3+2\leq n+1$.\\
\indent $\bullet$\; If $d_3+2=n+1$, then $n=7$. Note that
$\det(L[1,2,4|1,3,5])=\begin{vmatrix}
n-4& -1 & -1\\
0& -1&0\\
0&-1& x\end{vmatrix}=-x(n-4)$, where $x=-1$, or 0. Since
$(n-3)\mid-x(n-4)$ then $x=0$. So $v_4\not\sim v_5$. In the same
way, we can see  that $v_4\not\sim v_7$, $v_5\not\sim v_6$ and
$v_7\not\sim v_6$. So there is no edges between the vertices
$v_1,v_5,v_7$ and  $v_2, v_4, v_6$. Moreover, we have $d_5\leq n-4$.
Note that $\det(L[2,3,5|3,5,7])=\begin{vmatrix}
-1& 0& 0\\
-1& -1&-1\\
-1&d_5& y\end{vmatrix}=y-d_5$, where $y=-1$, or 0. From $(n-3)\mid
(y-d_5)$, we can get $d_5=n-4$ and $y=-1$. So $v_5\sim v_7$.
In the same way, we can get $v_4\sim v_6$. Thus $G=K_7-K_{3,3}$ (See Fig.1).\\
\indent $\bullet$\; If $d_3+2=n$ or $n-1$, then $n=6$ or 5 respectively, impossible.\\
\indent $\bullet$\; If $d_3+2=n-3$, then $d_3=n-5$. Then there
exists a vertex $v_8$ such that $v_1\sim v_8$, $v_2\sim v_8$ and
$v_3\not\sim v_8$. Thus we have
$\det(L[1,2,3|6,7,8])=\begin{vmatrix}
0&-1&-1\\
-1&0&-1\\
-1&-1&0
\end{vmatrix}=-2$. From $(n-3)|2$ we get $n=5$. Impossible.\\
\indent
Now we assume that $case$ 7 occurs. We know only its subcase 1 and subcase2 may  occur.
If its $subcase$ 1 occurs, then  $\det(L[1, 2,
3|1, 4, 5])=\left|
  \begin{array}{ccc}
    d_1 & -1 & -1 \\
    0 & -1 & -1 \\
    -1 & 0 & -1 \\
  \end{array}
\right|=d_1=n-2$. So $n-3$ divides $n-2$,  impossible.
If its $subcase$ 2 occurs, if we regard the
vertices $v_4$ as $v_1$, $v_5$ as  $v_2$, $v_1$ as $v_3$, $v_3$ as
$v_4$, then we are back in case 2. The required graphs have been determined.\\
\indent If $case$ 8 occurs, then only its  subcase 1 may occur.
Of course, we are back in the subcase 1 of case 7 and the required graphs have been determined.\\

\end{document}